\documentclass[leqno]{article}
\usepackage{amssymb}
\usepackage{amsmath, enumerate, vmargin}

\makeatletter
\renewcommand\section{\@startsection {section}{1}{\z@}%
 {-3.5ex \@plus -1ex \@minus -.2ex}%
 {2.3ex \@plus.2ex}%
 {\center \normalfont\large\bfseries}}
\makeatother

\setmarginsrb{3cm}{3cm}{3cm}{3cm}{1mm}{6mm}{0mm}{10mm}

\newtheorem{thm}{Theorem}

\newtheorem{cor}[thm]{Corollary}
\newtheorem{lem}[thm]{Lemma}
\newtheorem{defi}[thm]{Definition}
\newtheorem{remark}[thm]{Remark}
\newtheorem{example}[thm]{Example}
\newtheorem{pb}[thm]{Problem}

\newcommand{\real}{{\mathbb R}}
\newcommand{\nat}{{\mathbb N}}
\newcommand{\ent}{{\mathbb Z}}

\newcommand{\A}{{\mathcal A}}

\newcommand{\D}{{\mathcal D}}
\newcommand{\E}{{\mathcal E}}

\newcommand{\M}{{\mathcal M}}
\newcommand{\N}{{\mathcal N}}

\newcommand{\R}{{\mathcal R}}
\renewcommand{\S}{{\mathcal S}}
\newcommand{\Z}{{\mathcal Z}}

\newcommand{\f}{\varphi}
\newcommand{\s}{\sigma}

\newcommand{\n}{\noindent}
\newcommand{\pf}{\noindent{\it Proof.~~}}
\newcommand{\cqd}{\hfill$\Box$}
\newcommand{\be}{\begin{eqnarray*}}
\newcommand{\ee}{\end{eqnarray*}}
\newcommand{\beq}{\begin{equation}}
\newcommand{\eeq}{\end{equation}}

\begin{document}

%%%%%%%%%%%%%%%%%%%%%%%%%%%%%%%%%%%%%%%%%%%%%%
%%%%%%%  Title
%%%%%%%%%%%%%%%%%%%%%%%%%%%%%%%%%%%%%%%%%%%%%

\title{On the maximality of subdiagonal algebras}
\author{Quanhua Xu}
\date{}
\maketitle

%%%%%%%%%%%%%%%%%%%%%%%%%%%%%%%%%%%%%%%%%%%%%
%%%%%%%%%%%%%%%%%%%%%%%%%%%%%%%%%%%%%%%%%%%%%

\begin{abstract}
 We consider Arveson's problem on the
maximality of subdiagonal algebras. We prove that a subdiagonal
algebra is maximal if it is invariant under the modular group of a
faithful normal state which is preserved by the conditional
expectation associated with the subdiagonal algebra.
\end{abstract}

%%%%%%%%%%%%%%%%%%%%%%%%%%%%%%%%%%%%%%%%%%%
%%%%%%%%%%%%%%%%%%%%%%%%%%%%%%%%%%%%%%%%%%%
 %\setcounter{page}{0}
 %\setcounter{section}{-1}
 %\thispagestyle{empty}

%%%%%%%%%%%%%%%%%%%%%%%%%%%%%%%%%%%%%%%%%%%
%%%%%%%%%%%%%%%%%%%%%%%%%%%%%%%%%%%%%%%%%%%
\makeatletter
\renewcommand{\@makefntext}[1]{#1}
\makeatother \footnotetext{\noindent Laboratoire de
Math\'{e}matiques,
 Universit\'{e} de Franche-Comt\'{e},
 25030 Besan\c con, cedex - France\\
 qx@math.univ-fcomte.fr\\
 2000 {\it Mathematics subject classification:}
  46L10,  47D25\\
 {\it Key words and phrases}:
  Subdiagonal algebra, maximality, modular
group, conditional expectation}

%%%%%%%%%%%%%%%%%%%%%%%%%%%%%%%%%%%%%%%%%%%
%%%%%%%%%%%%%%%%%%%%%%%%%%%%%%%%%%%%%%%%%%%
\section{Introduction}

Let ${\M}$ be a von Neumann algebra. Let ${\E}$ be a normal
faithful conditional expectation from ${\M}$ onto a von Neumann
subalgebra ${\D}$ of ${\M}.$ A $\sigma$-weakly closed subalgebra
${\A}$ of ${\M}$ is called a {\it subdiagonal algebra} in ${\M}$
with respect to ${\E}$ if the following conditions are satisfied
 \begin{enumerate}[(i)]
 \item  ${\A} + {\A}^{\ast}$ is $\sigma$-weakly dense in
${\M}$;
 \item  ${\E}$ is multiplicative on ${\A}$;
 \item  ${\A} \cap {\A}^{\ast} = {\D}$, where ${\A}^{\ast} =
\{ x^{\ast}\ :\ x \in {\A} \}.$
\end{enumerate}
\noindent ${\D}$ is then called the {\it diagonal} of ${\A}.$

\smallskip

This notion was introduced by Arveson in \cite{ar} with the
perspective to give  a  unified theory of non-selfadjoint operator
algebras, including the algebra of bounded analytic matrix valued
(or more generally, operator valued) functions and nest algebras.
One fundamental result proved in \cite{ar} is an inner-outer type
factorization, which extends significantly the previous
inner-outer factorization for analytic matrix valued functions
obtained independently by Helson - Lowdenslager \cite{helo} and
Wiener - Masani \cite{wm} (see also \cite{he}). This theorem was
further generalized and studied in many related contexts (see
\cite{js} and for more references therein). On the other hand, as
the well-known Sz\"ego inner-outer factorization in the theory of
the classical Hardy spaces, this factorization is central for the
development of the non-commutative Hardy space theory (cf.
\cite{mar, marw, marw2}).

In all these works, and in fact since the creation of the theory
of subdiagonal algebras by Averson, a certain maximality
assumption has always played a preeminent r\^ole. Recall that a
subdiagonal algebra ${\A}$ with respect to ${\E}$ is said to be
{\it maximal} if ${\A}$ is properly contained in no larger
subdiagonal algebra with respect to ${\E}.$ It was proved in
\cite{ar} that any subdiagonal algebra ${\A}$ is contained in a
unique maximal subdiagonal algebra, denoted by ${\A}_{\max}$,
which is described by
 $${\A}_{\max} = \{ x \in {\M}\ :\ {\E} ( a x b) = 0,\
 \forall\; a \in {\A},\ \forall\; b \in {\A}_{0} \},$$
 where
 $${\A}_{0} = \{ a \in {\A}\ :\ {\E} (a) = 0 \}.$$
Many known examples of subdiagonal algebras are maximal. A long
standing open problem raised by Arveson in \cite{ar} is that
whether every subdiagonal algebra is \underbar{automatically}
maximal.

Only more than two decades later that Exel \cite{ex} gave a
partial solution for this problem: if there is a normal faithful
tracial state $\tau$ on ${\M}$ such that $\tau \circ {\E} = \tau$
(in this case ${\A}$ is called a {\it finite subdiagonal
algebra}), then ${\A}$ is maximal. In fact, Exel's arguments show
a little bit more, namely, that every subdiagonal algebra of a
finite von Neumann algebra is automatically maximal. By the way,
we recall another problem posed in \cite{ar}, still unsolved too,
is that whether a subdiagonal algebra of a finite and
$\sigma$-finite von Neumann algebra is a finite subdiagonal
algebra.

Very recently, Ji, Ohwada and Saito proved in \cite{jos} that if
${\A}$ is a maximal subdiagonal algebra in a $\sigma$-finite von
Neumann algebra ${\M}$ with respect to ${\E},$ then ${\A}$ is
invariant under the modular automorphism group
$\sigma_{t}^{\varphi}$ of every ${\E}$-invariant normal faithful
state $\varphi$ on ${\M}.$ Recall that $\varphi$ is
${\E}$-invariant if $\varphi \circ {\E} = \varphi.$ They then
asked that whether the converse is true. Let us explicitly state
this question as follows (see \cite[Question 2.7]{jos}).

\smallskip\n{\bf Question.} {\em Let ${\A}$ be a subdiagonal
algebra of a $\sigma$-finite von Neumann algebra ${\M}$ with
respect to ${\E}.$ Assume that ${\A}$ is
$\sigma^{\varphi}_{t}$-invariant (i.e.,
$\sigma_{t}^{\varphi}({\A}) \subset {\A},\, \forall\;t \in
{\real})$ for every ${\E}$-invariant normal faithful state
$\varphi$ on ${\M}.$ Is ${\A}$ maximal}?

\smallskip\n The aim of this note is to answer this question in the
affirmative. Below is our main result.

\begin{thm}\label{thm}
  Let ${\M}$ be a $\sigma$-finite von Neumann algebra and
${\E}$ a normal faithful conditional expectation from ${\M}$ onto
a von Neumann subalgebra ${\D} \subset {\M}.$ Let ${\A}$ be a
subdiagonal algebra of ${\M}$ with respect to ${\E}.$ If there is
a normal faithful state $\varphi$ on ${\M}$ such that $\E$
commutes with $\sigma_{t}^{\varphi}$ (i.e.,
 $\sigma_{t}^{\varphi}\circ\E
 ={\E}\circ\sigma_{t}^{\varphi}$ for all $t\in\real$) and
 $\A$ is $\sigma_{t}^{\varphi}$-invariant,
 then ${\A}$ is maximal.
 \end{thm}

\n{\bf Remark.} It is classical that if $\f$ is $\E$-invariant,
then $\E$ and $\sigma_{t}^{\varphi}$  commute (\cite[1.4.3]{co}).

\smallskip

The remainder of the note is essentially devoted to the proof of
the theorem above. Our strategy  is to reduce the present
situation to that of finite von Neumann algebras, and then to use
Exel's theorem quoted previously. The key ingredient of this
reduction is an unpublished important result of Haagerup. It
roughly says that every von Neumann algebra can be embedded, in an
appropriate way, into a large von Neumann algebra, which is a kind
of inductive limit of some nice finite von Neumann subalgebras. In
the next section, we will recall this reduction theorem of
Haagerup and the construction of these nicely disposed
subalgebras. The proof of the above theorem will be given in
section 3. Section 4 contains a generalization to weights instead
of states.

%%%%%%%%%%%%%%%%%%%%%%%%%%%%%%%%%%%%%%%%%%%
%%%%%%%%%%%%%%%%%%%%%%%%%%%%%%%%%%%%%%%%%%%

\section{ Haagerup's reduction theorem}

In this section we recall an important unpublished theorem due to
Haagerup \cite{haag-red}. It states that any von Neumann algebra
can be embedded, as the image of a normal faithful conditional
expectation, into a large von Neumann algebra which is generated
by an increasing family of finite subalgebras, each of which is
the image of a normal conditional expectation. Haagerup's original
intention is to approximate his non-commutative $L^{p}$-spaces
based on type III von Neumann algebras by those constructed from a
trace. This approximation theorem on Haagerup non-commutative
$L^{p}$-spaces is very important in non-commutative analysis. In
many situations, it permits to consider only non-commutative
$L^{p}$-spaces associated with traces. We refer to \cite{jx} for
more recent applications of Haagerup's reduction theorem to
non-commutative martingale and ergodic theories. Note that
\cite{jx} also contains a reproduction of Haagerup's unpublished
manuscript  \cite{haag-red}.

The main tool of Haagerup's construction is crossed products. Our
references for crossed products are \cite{kar, stra}. Throughout,
$G$ will denote the discrete subgroup $\bigcup_{n \geq 1} 2^{-n}
{\ent}$ of ${\real}.$ Let ${\M}$ be a von Neumann algebra acting
on a Hilbert space $H$ and $\varphi$ a normal faithful state on
${\M}.$ We consider the crossed product
 ${\M}\rtimes_{\sigma^{\varphi}}G$ of ${\M}$ by $G$
 with respect to $\sigma^{\varphi}.$ In the sequel, we will denote this crossed
product by ${\R}.$ Recall that ${\R}$ is a von Neumann algebra on
$\ell^{2} (G, H)$ generated by the operators $\pi (x), x \in {\M}$
and $\lambda(t),\ t \in G,$ which are defined by
 $$\big(\pi(x) \xi\big)(s)=\sigma^{\varphi}_{-s}(x)\xi(s),\quad
 \big(\lambda (t)\xi\big)(s)=\xi (s-t),\ \
 s \in G,\ \xi \in \ell^{2}(G, H).$$
Note that $\pi$ is a normal faithful representation of ${\M}$ on
$\ell^{2}(G, H).$ Thus we  will  identify $\pi({\M})$ and ${\M}$
whenever possible. The operators $\pi(x)$ and $\lambda(t)$ satisfy
the following commutation relation:
 \beq
 \lambda(t)\pi(x)\lambda(t)^{\ast} = \pi (\sigma^\f_{t}(x)),\ \
  t\in G, \ x \in {\M}.
  \eeq
Let $\widehat \varphi$ be the dual weight of $\varphi$ on ${\R}.$
Then $\widehat \varphi$ is again a normal faithful state on ${\R}$
uniquely determined by
 \beq
 \widehat \varphi(\lambda (t) x) =\left\{ \begin{array}{ll}
 \varphi(x)  & \textrm{ if}\ t = 0\\
 0 & \textrm{ otherwise}
 \end{array}\right.,
\quad  x \in {\M},\ t \in G.
 \eeq
In particular, $\widehat \varphi \big|_{\M} = \varphi.$ The
modular automorphism group of $\widehat \varphi$ is uniquely
determined by
 \beq
 \sigma_{t}^{\widehat \varphi} (x) =
 \sigma_{t}^{\varphi}(x),\quad
 \sigma_{t}^{\widehat \varphi} (\lambda (s)) = \lambda(s),\ \
 x \in {\M},\ t, s \in G.
 \eeq
Consequently, $\sigma_{t}^{\widehat \varphi}\big|_{{\M}} =
\sigma_{t}^{\varphi},$ and so $\sigma_{t}^{\widehat \varphi}({\M})
= {\M}$ for all $t \in {\real}.$ It also follows that
 \beq
 \sigma_{t}^{\widehat \varphi}(x)
 = \lambda(t) x\lambda(t)^{\ast},\ \
 x \in {\R},\ t \in G.
 \eeq
 It is classical that there is a unique normal faithful conditional
expectation $\Phi$ from ${\R}$ onto ${\M}$ determined by
 \beq
 \Phi\big(\lambda(t) x\big) = \left\{\begin{array}{ll}
 x  & \textrm{ if }\ t = 0\\
  0 & \textrm{ otherwise}
  \end{array}\right.,
  \quad x \in {\M},\ t \in G.
 \eeq
 By (2), (3) and
(5), we deduce that
 \beq
 \widehat \varphi \circ \Phi = \widehat \varphi \quad
 \hbox{and}\quad
 \sigma_{t}^{\widehat \varphi} \circ \Phi = \Phi \circ
 \sigma_{t}^{\widehat \varphi},\ \ t \in {\real}.
 \eeq
With these notations, Haagerup's reduction theorem asserts that
there is an increasing sequence $({\R}_{n})_{n \geq 1}$ of von
Neumann subalgebras of ${\R}$ with the following properties:
 \begin{enumerate}[(i)]
 \item  each ${\cal R}_{n}$ is finite ;
 \item  $\bigcup_{n\geq 1}\ {\R}_{n}$ is $\s$-weakly dense in ${\R}$ ;
 \item for every $n \geq 1$ there is a normal faithful conditional
 expectation $\Phi_{n}$ from ${\R}$ onto ${\R}_{n}$ such that
 \beq
 \widehat \varphi \circ \Phi_{n} = \widehat \varphi,\quad
 \sigma_{t}^{\widehat \varphi} \circ \Phi_{n} = \Phi_{n} \circ
 \sigma_{t}^{\widehat \varphi},\quad
 \Phi_{n} \circ \Phi_{n+1} =
 \Phi_{n},\ \ n \geq 1,\ t \in {\real}.
 \eeq
 \end{enumerate}

Note that a normal conditional expectation satisfying the first
equality in (7) is unique. Since $\Phi_n \circ \Phi_{n+1}$ is also
conditional expectation under which $\widehat \varphi$ is
invariant, this uniqueness implies $\Phi_{n} \circ \Phi_{n+1} =
\Phi_{n},$ that is, the third equality in (7) is a consequence of
the first. Note that the second equality is also a consequence of
the first by Connes' classical result already quoted before.

In Haagerup's construction, ${\R}_{n}$ is the centralizer of a
normal faithful state $\varphi_{n}$ on ${\R}$ such that its
modular automorphism group $\sigma_{t}^{\varphi_{n}}$ is periodic
of period $2^{-n}.$ In the sequel, we will need the precise form
of $\varphi_{n}.$ Thus let us briefly recall this construction.

For a von Neumann algebra ${\N}$ and a normal faithful state
$\psi$ on ${\N}$ we denote, as usual, by ${\Z} ({\N})$ the center
of ${\N}$ and by ${\N}_{\psi}$ the centralizer of $\psi$ in
${\N}.$ Recall that ${\N}_{\psi}$  is the algebra of the fixed
points of $\sigma^{\psi}_{t}$. By (4), $\lambda(t) \in
{\Z}({\R}_{\widehat \varphi})$ for all $t \in G.$ For any given $n
\in {\nat},$ by functional calculus, there is $b_{n} \in {\Z}
({\R}_{\widehat \varphi})$ such that
 $$ 0 \leq b_{n} \leq 2 \pi \quad \mbox{ and } \quad
 e^{i b_{n}} = \lambda (2^{-n}).$$
Set $a_{n} = 2^{n} b_{n}.$ Then again $a_{n} \in
{\Z}({\R}_{\widehat \varphi}),\ n \geq 1.$ The desired state
$\varphi_{n}$ is defined as
 \beq
 \varphi_{n}(x) = {1\over \widehat\varphi(e^{-a_{n}})}
\; \widehat \varphi(e^{-a_{n}} x)
  ,\ \
 x \in {\R},\ n \geq 1.
 \eeq
Since $a_{n} \in {\R}_{\widehat \varphi},$
 \beq\sigma_{t}^{\varphi_{n}} (x) = e^{-it a_{n}}
 \sigma_{t}^{\widehat \varphi}(x) e^{it a_{n}},\ \
 x \in {\cal R},\ t\in {\real},\ n \geq 1.
 \eeq
Then by (4) and the definition of $a_n$,
$\sigma_{t}^{\varphi_{n}}$ is $2^{-n}$-periodic. Let ${\R}_n =
{\R}_{\varphi_{n}}.$ Then $\varphi_{n} \big|_{{\R}_{n}}$ is a
normal faithful tracial state on ${\R}_{n}$, and so ${\R}_{n}$ is
a finite von Neumann subalgebra of ${\R}.$

Define $\Phi_{n}\ :\ {\R} \to {\R}_{n}$ by
 $$\Phi_{n}(x) = 2^{n}\int^{2^{-n}}_{0}\sigma_{t}^{\varphi_{n}}
 (x)dt,\quad x \in {\R}.$$
By the $2^{-n}$-periodicity of $\sigma_{t}^{\varphi_{n}},$ we have
 \beq
 \Phi_{n}(x) = \int^{1}_{0} \sigma_{t}^{\varphi_{n}} (x)dt,
 \quad x \in {\R}.
 \eeq
Then it is routine to check that $\Phi_{n}$ is a normal faithful
conditional expectation satisfying (7). Hence to prove Haagerup's
reduction theorem mentioned above it remains to show that
$({\R}_{n})$ is increasing and the union of the ${\R}_{n}'s$ is
$\s$-weakly dense in ${\R}.$ We refer the reader to
\cite{haag-red, jx} for more details.

%%%%%%%%%%%%%%%%%%%%%%%%%%%%%%%%%%%%%%%%%%%
%%%%%%%%%%%%%%%%%%%%%%%%%%%%%%%%%%%%%%%%%%%

\section{ The proof}

This section is devoted to the proof of Theorem \ref{thm}.
Throughout this section, ${\M}, {\D}, {\E}, {\A}$ and $\varphi$
will be fixed as in that theorem.  ${\R}$ will be the crossed
product ${\M} \rtimes_{\sigma^{\varphi}}G$ as in the last section,
and we will keep all notations introduced there. The idea of the
proof is to first lift $\A$ to a subdiagonal algebra in $\R$, then
compress the latter  to a subdiagonal algebra in $\R_n$ by the
conditional expectation $\Phi_n$, and finally come back to $\A$ by
passing to limit as $n\to\infty$.

For easy later reference let us state the commutation assumption
on $\E$ and $\sigma_{t}^{\varphi}$ as follows
 \beq
 \sigma_{t}^{\varphi} \circ {\E} = {\E} \circ
  \sigma^{\varphi}_{t}, \ \ t\in\real.
  \eeq
This implies that  ${\D}$ is $\sigma_{t}^{\varphi}$-invariant and
$\sigma^{\varphi}_{t} \big|_{{\D}}$ is exactly the modular
automorphism group of $\varphi \big|_{{\D}}.$ Consequently, we do
not need to distinguish $\varphi$ and $\varphi\big| _{{\D}},$
$\sigma_{t}^{\varphi}$ and $\sigma_{t}^{\varphi}\big|_{{\D}}$,
respectively.  Now let ${\S} = {\D} \rtimes_{\sigma^{\varphi}} G.$
Then ${\S}$ is naturally identified as a von Neumann subalgebra of
${\R},$ generated by all operators $\pi(x),\, x \in {\D}$ and
$\lambda(t),\, t \in G.$ The dual weight of $\varphi \big|_{{\D}}$
on ${S}$ is equal to $\widehat \varphi \big|_{{\S}}.$ Again, we
will denote this restriction by the same symbol $\widehat
\varphi.$ It is not hard to extend ${\E}$ to a normal faithful
conditional expectation $\widehat {\E}$ from ${\R}$ onto ${\S},$
which is uniquely determined by
 \beq
 \widehat {\E} (\lambda(t) x)
 = \lambda (t)\, {\E}(x),\quad
 x \in {\M},\quad t \in G.
 \eeq
The reader is referred to \cite{jx} for details and for more
extensions of this type. By (4), (11) and (12), we deduce
 \beq
 \sigma_{t}^{\widehat \varphi} \circ \widehat {\E}
 = \widehat{\E} \circ \sigma_{t}^{\widehat \varphi},
 \quad t \in G.
 \eeq
On the other hand,  using (9), (13) and the fact that $a_{n} \in
{\S}$ and $\widehat {\E}$ is a conditional expectation with
respect to ${\S}$, we get
 \beq
 \sigma_{t}^{\varphi_{n}} \circ \widehat {\E}
 = \widehat {\E} \circ \sigma_{t}^{\varphi_{n}},\quad
  t \in {\real},\ n \geq 1.
  \eeq
Hence by the definition (10) of the conditional expectation
$\Phi_{n}\ :\ {\R} \to {\R}_{n},$ we deduce
 \beq
 \Phi_{n} \circ \widehat {\E}
 =\widehat{\E}\circ\Phi_{n},\quad n\geq 1.
 \eeq
In particular, ${\R}_{n}$ and ${\S}$ are respectively $\widehat
{\E}$-invariant and $\Phi_n$-invariant.

Now let ${\S}_{n} = {\S}_{\varphi_{n}},\ n \geq 1.$ Then clearly,
${\S}_{n} = {\R}_{n} \cap {\S}$ for every $n\geq 1$. Also note
that $\Phi_{n} \big|_{{\S}}$ and $\widehat {\E} \big|_{{\R}_{n}}$
are normal faithful conditional expectations from ${\S}$ onto
${\S}_{n}$, respectively, from ${\R}_{n}$ onto ${\S}_{n}$.
$({\S}_{n})_{n \geq 1}$ and $(\Phi_{n}\big|_{{\S}})_{n \geq 1}$
are the increasing sequences of von Neumann subalgebras of ${\S}$
and respectively the sequence of the corresponding conditional
expectations given by Haagerup's construction presented in the
last section relative to $({\D},\, \varphi \big|_{{\D}})$ instead
of $({\M}, \varphi).$
 Again, we will denote these
restriction mappings by the same symbols as the mappings
their-selves when no confusion can occur.

Since ${\A}$ is $\sigma^{\varphi}_{t}$-invariant, by (1), the
family of all linear combinations on $\lambda(t)\, \pi(x),$ $t \in
G,$ $x \in {\A},$ is a $\ast$-subalgebra of ${\R}.$ Let $\widehat
{\A}$ be its $\s$-weakly closure in ${\R}$ and ${\A}_{n} =
\widehat {\A} \cap {\R}_{n}$. The following lemmas show that
$\widehat {\A}$ (resp. ${\A}_{n}$) is a subdiagonal algebra with
respect to $\widehat {\E}$ (resp. $\widehat {\E}
\big|_{{\R}_{n}}$).

\begin{lem}\label{lem1}
 $\widehat {\A}$ is a subdiagonal algebra of ${\R}$ with
respect to $\widehat {\E}$.
\end{lem}

\pf We first prove that $\widehat {\A} + \widehat {\A}^{\ast}$ is
$\s$-weakly dense in ${\R}.$ For this it suffices to show that for
any $t \in G$ and $x \in {\M},$ $\lambda(t)\, \pi(x)$ is the limit
of elements in $\widehat {\A} + \widehat {\A}^{\ast}.$ Since ${\A}
+ {\A}^{\ast}$ is $\s$-weakly dense in ${\M},$ there are $a_{i},
b_{i} \in {\A}$ such that
 $$x = \lim_{i}(a_{i} + b_{i}^{\ast})\quad \sigma{\rm -weakly}.$$
Since $\pi$ is normal,
 $$\pi(x) = \lim_{i} (\pi(a_{i}) + \pi (b_{i})^{\ast})\quad
 \sigma{\rm -weakly}.$$
Therefore,
 $$\lambda(t)\pi(x) = \lim_{i} (\lambda(t)\pi(a_{i}) +
 \lambda(t)\pi(b_{i})^{\ast}) \quad \sigma{\rm -weakly}.$$
This is the desired limit.

Next we show that $\widehat {\E}$ is multiplicative on $\widehat
{\A}.$ To this end we note that by (12), for any $s, t \in G$ and
$x, y\in {\A}$
 \be
 \widehat {\E}\big(\lambda(s)\pi(x)\pi(y)\lambda(t)\big)
 &=&\lambda(s)\pi\big({\cal E}(x y)\big)\lambda(t)\\
 &=&\lambda(s)\pi\big({\cal E}(x){\cal E}(y)\big)\lambda(t)\\
 &=& \widehat {\cal E}\big(\lambda(s)\pi(x)\big)
 \widehat{\cal E}\big(\pi(y)\lambda(t)\big),
 \ee
 where we have used the mutliplicativity of ${\E}$ on
${\A}.$ Then the linearity and normality of $\widehat {\E}$ imply
the mutliplicativity of $\widehat {\A}$ on $\widehat {\A}.$

Thus it remains to show $\widehat {\A} \cap \widehat {\A}^{\ast} =
{\S}.$ To this end, we will use the matrix representation
$(x_{s,t})_{s,t \in G}$ of an element $x \in B(\ell^{2}(G, H))$ in
the natural basis of $\ell^{2}(G).$ It is well-known that $x \in
{\R}$ iff there is a function $X\ :\ G \to {\M}$ such that
 $$x_{s, t} = \sigma^{\varphi}_{-s}\big(X( s\, t^{-1})\big),
 \quad s, t \in G.$$
 (cf. \cite[section 22.1]{stra}). Clearly, this
function $X$ is unique. Now we claim that if $x \in \widehat
{\A},$ then $X(t) \in {\A}$ for all $t \in G.$ Indeed, this is
clear if $x = \lambda (t_{0})\, \pi(x_{0})$ for some $t_{0} \in G$
and $x_{0} \in {\A}.$ It then follows that the claim is true if
$x$ is a linear combination of $\lambda(t)\, \pi(y),\ t \in G,\ y
\in {\A}.$ For a general $x \in \widehat {\A},$ there is a net $\{
x_{i} \}$ of linear combinations on $\lambda(t)\, \pi(y),\, t \in
G,\, y \in {\A},$ such that
 $$x = \lim_{i}x_{i} \hskip 0.5cm \sigma{\rm -weakly}.$$
If $X_{i}$ denotes the function corresponding to $x_{i},$ then
clearly
 $$X(t) = \lim_{i}X_{i}(t) \hskip 0.5cm \sigma{\rm -weakly},\ t \in
 G.$$
Hence by the $\sigma$-weak closedness of ${\A},$ we conclude that
$X(t) \in \A$ for all $t \in G,$ proving our claim.

Similarly, if $x \in \widehat {\A}^{\ast},$ then $X(t) \in
{\A}^{\ast}$ for all $t \in G.$ Now let $x \in \widehat {\A} \cap
\widehat {\A}^{\ast}.$ Then $X(t) \in {\A} \cap {\A}^{\ast} =
{\D}$ for all $t \in G.$ Therefore, $x \in {\S},$ and so $\widehat
{\A} \cap \widehat {\A}^{\ast} \subset {\S}.$ The converse
inclusion is trivial. Thus $\widehat {\A} \cap \widehat
{\A}^{\ast} = {\cal S}$. Therefore $\widehat {\A}$ is a
subdiagonal algebra with respect to $\widehat {\E}.$ \cqd

\begin{lem}\label{lem2}
 Every ${\A}_{n}$ is a finite subdiagonal algebra in
$\R_n$ with respect to $\widehat {\E} \big|_{{\R}_{n}}.$
\end{lem}

\pf  Since $\widehat {\E}$ is multiplicative on $\widehat {\A},$
$\widehat {\E} \big|_{{\R}_{n}}$ is multiplicative on ${\A}_{n}.$
On the other hand,
 $${\A}_{n} \cap {\A}_{n}^{\ast}
 = \widehat {\A} \cap\widehat {\A}^{\ast}
 \cap {\R}_{n} = {\S} \cap {\R}_{n} = {\S}_{n}.$$
Thus it remains to show the $\s$-weak density of ${\A}_{n} +
{\A}^{\ast}_{n}$ in ${\R}_{n}.$ Let $x \in {\R}_{n}.$ Since
$\widehat {\A} + \widehat {\A}^{\ast}$ is $\s$-weakly dense in
${\R},$ there are $a_{i}, b_{i} \in \widehat {\A}$ such that
 $$x = \lim_{i}(a_{i} + b_{i}^{\ast}) \hskip 0.5cm
 \sigma{\rm-weakly}.$$
Then by the normality of $\Phi_{n},$ we have
 $$x = \Phi_{n}(x) = \lim_{i}(\Phi_{n}(a_{i}) +
 \Phi_{n}(b_{i})^{\ast}) \ \ \sigma{\rm -weakly}.$$
However, by (9), (10) and the assumption that ${\A}$ is
$\sigma_{t}^{\varphi}$-invariant, we easily deduce that $\widehat
{\A}$ is $\Phi_{n}$-invariant for all $n \geq 1.$ Hence,
$\Phi_{n}(a_{i}),\ \Phi_{n}(b_{i}) \in \widehat {\A} \cap {\cal
R}_{n} = {\A}_{n}.$ It follows that ${\A}_{n} + {\A}^{\ast}_{n}$
is $\sigma$-weakly dense in ${\R}_{n}.$ Thus ${\A}_{n}$ is a
subdiagonal algebra with respect to $\widehat {\E}
\big|_{{\R}_{n}}$. Note that as a by-product we have also proved
${\A}_{n}=\Phi_n(\R_n)$.  \cqd

\smallskip

We recall that if ${\A}$ is a subdiagonal algebra in ${\M}$ with
respect to ${\E},$ then the maximal subdiagonal algebra containing
${\A}$ is
 $${\A}_{\max} = \{ x \in {\M}\ :\
 {\E}({\cal A}x{\A}_{0}) = {\E}({\A}_{0}x
 {\A}) =0\}.$$

\begin{lem}\label{lem3}
 $\widehat {\A}$ is maximal.
 \end{lem}

\pf We must show $(\widehat {\A})_{\max} = \widehat {\A}.$ Let $x
\in (\widehat {\A})_{\max}.$ Set $x_{n} = \Phi_{n}(x), n \geq 1.$
We claim that $x_{n} \in ({\A}_{n})_{\max}.$ Indeed, let $a, b \in
{\A}_{n}$ with $\widehat {\E} (b) = 0.$ Then $a, b \in {\A} \cap
{\R}_{n}.$ Since $\Phi_{n}$ is a conditional expectation with
respect to ${\R}_{n},$ by (15), we have
 $$\widehat {\E} (ax_{n} b)
 = \widehat {\E}(a\Phi_{n}(x)b)
 = \widehat {\E} (\Phi_n (axb))
 = \Phi_n (\widehat {\E} (axb)) = 0.$$
 This yields our claim. However, by Lemma \ref{lem2} and Exel's
theorem, ${\A}_{n}$ is maximal. Hence $x_{n} \in \widehat {\A}$
for $n \geq 1.$ On the other hand, (7) implies that $x_{n} \to x$
$\sigma$-weakly. Since $\widehat {\A}$ is $\sigma$-weakly closed,
we conclude that $x \in \widehat{\A}.$ Therefore, $\widehat {\A}$
is maximal. \cqd

\smallskip

Finally, we are ready to prove our main theorem.

\smallskip

\n{\it Proof of Theorem \ref{thm}.}  Applying  the preceding
discussion to ${\A}_{\max}$ in the place of ${\A}$,  we get a
subdiagonal algebra $\widehat{{\A}_{\max}}$ of ${\R}$ with respect
to $\widehat{\E}.$ Since ${\A} \subset {\A}_{\max}$, $\widehat
{\A} \subset \widehat{{\A}_{\max}}$. However, by Lemma \ref{lem3},
$\widehat {\A}$ is maximal. Hence $\widehat {\A} =
\widehat{{\A}_{\max}}$. Consequently, for any $x \in {\A}_{\max},\
\pi(x) \in \widehat{{\A}_{\max}}= \widehat {\A}.$ Then
necessarily, $x \in {\A}.$ Thus ${\A} = {\A}_{\max},$ and so
${\A}$ is maximal.\cqd

%%%%%%%%%%%%%%%%%%%%%%%%%%%%%%%%%%%%%%%%%%%
%%%%%%%%%%%%%%%%%%%%%%%%%%%%%%%%%%%%%%%%%%%

\section{ A generalization}

 It is not clear to the author at the time of this writing
whether the state $\f$ in Theorem \ref{thm} can be replaced by a
semifinite normal faithful weight (keeping all other assumptions).
The author is able to prove this only for normal faithful weights
whose restrictions to $\D$ are strictly semifinite. Recall that a
weight $\f$ on $\M$ is said to be strictly semifinite if there is
a family $\{\psi_j\}_{j\in J}$ of normal positive functionals
whose supports are  pairwise disjoint and such that
 $$\f=\sum_{j\in J}\psi_j.$$
This is equivalent to saying that $\f$ is semifinite on the
centralizer $\M_\f$. Our main theorem can be extended to weights
as follows.

 \begin{thm}\label{thmbis}
 Let ${\M}$ be a von Neumann algebra and ${\E}$ a
normal faithful conditional expectation from ${\M}$ onto a von
Neumann subalgebra ${\D} \subset {\M}.$ Let ${\A}$ be a
subdiagonal algebra of ${\M}$ with respect to ${\E}.$ If there is
a normal faithful weight $\varphi$ on ${\M}$ such that
$\f\big|_{\D}$ is strictly semifinite on $\D$, $\E$ commutes with
$\sigma_{t}^{\varphi}$ and
 $\A$ is $\sigma_{t}^{\varphi}$-invariant,
 then ${\A}$ is maximal.
 \end{thm}

As a corollary, we get the following generalization of Exel's
theorem to the semifinite case. See \cite{ji} for a related
result.

\begin{cor}
 Let ${\A}$ be a subdiagonal algebra of ${\M}$ with
respect to ${\E}.$ If there is a normal semifinite faithful trace
$\tau$ on ${\M}$ such that $\tau$ is semifinite on $\D$, then
${\A}$ is maximal.
 \end{cor}

The proof of Theorem \ref{thmbis} above can be reduced to the
state case via a standard way. Indeed, let $\f$ be a weight as in
the theorem and consider again the crossed product
$\R={\M}\rtimes_{\sigma^{\varphi}}G$. Using the strict
semifiniteness and the construction in section 2, one can prove
that there is an increasing family $\{ {\R}_{i} \}_{i\in I}$ of
$w^{\ast}$-closed $\ast$-subalgebras of $\R$ satisfying the
following properties :
 \begin{enumerate}[(i)]
 \item  each ${\R}_{i}$ is finite and $\sigma$-finite ;
 \item the union of all ${\R}_{i}$ is $w^{\ast}$-dense
 in ${\R}$ ;
 \item  the identity $p_{i}$ of ${\R}_{i}$ belongs to ${\R}_{\widehat\varphi}$ ;
 \item there is a normal conditional expectation $\Phi_{i}$
 from ${\R}$ onto ${\R}_{i}$ such that
 $$\widehat\varphi \circ \Phi_{i} = p_i \widehat\varphi
 p_{i} \quad \mbox{ and }\quad
 \sigma_{t}^{\widehat \varphi} \circ \Phi_{i}
 = \Phi_{i} \circ \sigma_{t}^{\widehat \varphi},\
 t \in {\real},\ i\in I.$$
 \item for all $i,j \in I$ with $i\leq j,$
 $$\Phi_{i} \circ \Phi_{j} = \Phi_{j} \circ
 \Phi_{i} = \Phi_{i}.$$
 \end{enumerate}

\ We refer to \cite{jx} for more details. Then repeating the
arguments in section 3, we can prove Theorem \ref{thmbis}. We omit
all details.

 \end{document}